\documentclass{amsart}

\usepackage[utf8]{inputenc}
\usepackage[T1]{fontenc}
\usepackage{textcomp}
\usepackage{amsfonts}
\usepackage{amsthm}
\usepackage{amssymb}
\usepackage{pst-node}
\usepackage[all]{xy}
\usepackage{mathtools}
\usepackage{chngcntr}
\usepackage{thmtools}
\usepackage{mathrsfs}
\usepackage{tcolorbox}
\usepackage{listings}
\usepackage{imakeidx}
\usepackage{hyperref}
\usepackage{color}
\usepackage{xcolor}
\usepackage{float}
\usepackage{tikz}
\usepackage{enumerate}
\usepackage[nameinlink, capitalize, noabbrev]{cleveref}
\usepackage{natbib}
\usetikzlibrary{fit}
\usepackage{tikz}
\usetikzlibrary{arrows.meta, positioning}

\hypersetup{
    colorlinks=true,
    linkcolor=teal,
    citecolor=magenta,
    }

\newtheorem{Theorem}{Theorem}

\newtheorem{proposition}{Proposition}[section]
\newtheorem{lemma}[proposition]{Lemma}

\newtheorem{theorem}[proposition]{Theorem}

\theoremstyle{definition}
\newtheorem{remark}[proposition]{Remark}

\def \<#1>{{\left\langle{#1}\right\rangle}}
\def\abs#1{\left\vert{#1}\right\vert}
\def\set#1{{\def\st{\;:\;}\left\{#1\right\}}}

\def\mat#1{\left(\begin{matrix}#1\end{matrix}\right)}

\def\NN{\mathbb N}
\def\ZZ{\mathbb Z}
\def\Z-{\overline{\mathbb Z}}
\def\Zl-{\overline{\mathbb Z}_\ell}

\def\QQ{\mathbb Q}
\def\Q-{\overline{\mathbb Q}}
\def\Ql-{\overline{{\mathbb Q}_\ell}}

\def\Fl-{\overline{{\mathbb F}_\ell}}
\def\FF{\mathbb F}

\def\PP{\mathbb P}

\DeclareMathOperator{\FPP}{FPP}

\DeclareMathOperator{\car}{char}

\DeclareMathOperator{\SL}{SL}

\DeclareMathOperator{\Stab}{Stab}

\DeclareMathOperator{\st}{st}
\DeclareMathOperator{\ev}{ev}

\DeclareMathOperator{\End}{End}
\DeclareMathOperator{\IIm}{Im}
\DeclareMathOperator{\id}{id}
\DeclareMathOperator{\Aut}{Aut}

\numberwithin{equation}{section}

\title[FPP of the geometric iterated Galois groups of Latt\`es maps]{Fixed-point proportion of the geometric iterated Galois groups of Latt\`es maps}
\author[J. Fari\~na-Asategui and S. Radi]{Jorge Fariña-Asategui and Santiago Radi}
\address{Jorge Fari\~{n}a-Asategui: Section de math\'{e}matiques, Universit\'{e} de Gen\`{e}ve, 
1211 Gen\`{e}ve 4, Switzerland}
 \email{Jorge.Farina-Asategui@unige.ch}
\address{Santiago Radi: Universidad de la República, 11200 Montevideo, Uruguay.
}
\email{santiradi19@gmail.com}
\thanks{The first author is supported by a postdoctoral grant in mathematics from the Knut och Alice Wallenberg foundation in Sweden. The second author is supported by a postdoctoral fellowship granted by PEDECIBA in Uruguay.}

\begin{document}

\begin{abstract}
Lattès maps, together with Chebyshev polynomials, are conjectured to be the only rational functions whose geometric iterated Galois groups have positive fixed-point proportion. Here, we explicitly compute the fixed-point proportion of the geometric iterated Galois group of a Lattès map, proving that it only depends on its degree and whether its postcritical orbits all intersect or not. In particular, we confirm that their fixed-point proportion is always positive.
\end{abstract}
\maketitle

\section{Introduction}
\label{section: introduction}

Given a rational function $f\in K(z)$ defined over a perfect algebraically closed field $K$, one obtains its \textit{geometric iterated Galois group} $G_\infty(K,f,t)$ as the Galois group corresponding to the infinite Galois extension $K(\bigcup_{n\ge 1}f^{-n}(t))/K(t)$ for a transcendental element $t$ over $K$. This group acts naturally on a regular rooted tree, which is no more than the tree of preimages of $t$ by the iterates of $f$. The \textit{fixed-point proportion} of $G_\infty(K,f,t)$ is the probability that a random element in $G_\infty(K,f,t)$, with respect to its unique Haar probability measure, fixes an infinite path in the tree.

Classifying the fixed-point proportion of geometric iterated Galois groups is a major problem in arithmetic dynamics, as it provides a very efficient tool to study different properties of iterates of polynomials over number fields and their residue fields; compare \cite{Bridy,JorgeSanti1, JorgeSanti2, JonesAMS,JonesComp,JonesLMS, JonesManes2012, Juul, Odoni1, Radi2026}. The recent breakthrough of the authors in \cite{JorgeSanti2} solved this problem completely in the polynomial case, proving that polynomials linearly conjugate to $\pm T_d$, where $T_d$ is the Chebyshev polynomial of degree $d$, are the only cases where the fixed-point proportion is positive. 

Based on the ideas of \cite{JorgeSanti2} and using the graph theoretic notion of transportation cost, the second author extended the results in \cite{JorgeSanti2} for rational functions, proving that the fixed-point proportion is zero for every rational function who is not a dynamical pullback and whose fixed-point process is a martingale; where a rational function is a \textit{dynamical pullback} if it is induced by the endomorphisms of a projective curve \cite{Radi2026inverseGalois}. This concept was introduced by Adams in \cite{Adams} in the study of projection-invariant subgroups of the geometric iterated Galois group. Such a curve, can either have genus zero or one by the Riemann-Hurwitz formula. In the latter case, this means that the rational function is induced by the endomorphisms of an elliptic curve and therefore it is a Lattès map. 

The main goal of this note is to show that Lattès maps yield positive fixed-point proportion. Our main result is the following:

\begin{Theorem}
\label{Theorem: fpp of Lattès}
Let $K$ be a perfect field and $f\in K(z)$ a separable Lattès map induced by an elliptic curve $E$. Then 
$$\mathrm{FPP}(G_\infty(K,f,t))>0.$$
Furthermore, when $f$ is tamely ramified and $E$ is not supersingular, the value $\mathrm{FPP}(G_\infty(K,f,t))$ can be calculated explicitly and only depends on the degree and the ramification portrait of $f$.
\end{Theorem}

Lattès maps are classified via their Thurston orbifolds (see \cref{section: classification} for further details), and they correspond to the Euclidean types:
$$(2,2,2,2),(2,4,4),(3,3,3),(2,3,6).$$ 

For Lattès maps with orbifold type other than (2,2,2,2), the value of the fixed-point proportion of their geometric iterated Galois groups only depends on the degree of the map, namely on whether it is divisible by $2,3$ or 6; see \textcolor{teal}{Theorems}~\ref{theorem: fpp (2,4,4)}, \ref{theorem: fpp (3,3,3)}, and \ref{theorem: fpp (2,3,6)}. For the orbifold type $(2,2,2,2)$, the fixed-point proportion also depends on whether the postcritical orbits intersect or not; see \cref{theorem: fpp (2 2 2 2)}.

\cref{Theorem: fpp of Lattès} is obtained by showing that for a Lattès map $f$, the fixed-point proportion is the sum of the reciprocals of finite indices obtained from the action of certain isogenies on the Tate module of the map inducing $f$. To show that the value of the fixed-point proportion only depends on the degree of the map and the topology of its postcritical orbits, we do a case-by-case analysis reducing the problem to studying the norm of certain ideals in an Euclidean domain associated to the elliptic curve.

\subsection*{\textit{\textmd{Organization}}}
In \cref{section: preliminaries}, we introduce basic notions of elliptic curves and Lattès maps. In \cref{section: general formula} we prove a general formula for the fixed-point proportion of the geometric iterated Galois group of a Lattès map; namely \cref{theorem: formula for FPP Lattès}. This already shows positivity of the fixed-point proportion for all Lattès maps. Lastly, in \cref{section: classification} we do a case-by-case analysis of each Lattès map and compute the corresponding fixed-point proportion.

\subsection*{AI Statement}
No AI was used neither in the writing nor in the proof of anything but the surjectivity claim in \cref{theorem: equivalence of End(E)-sets}. The model used was Claude Sonnet 5.5. The argument was simplified by the authors and it is only needed when the Lattès map is not defined over the complex numbers.

\subsection*{Remark}
Right before uploading this note, we came to know that Li had also proved positivity of the fixed-point proportion of Lattès maps in fields of characteristic zero in \cite{Li2026Lattes}. Our result was obtained independently of \cite{Li2026Lattes}. 

\section{Preliminaries on elliptic curves}
\label{section: preliminaries}

We first recall some results of elliptic curves that will be needed in the subsequent sections.

\subsection{Elliptic curves}

Let $E/K$ be an elliptic curve over a perfect algebraically closed field $K$. We denote $E(K)$ the set of points of $E$ and $\End(E)$ the set of isogenies of $E$, namely, endomorphisms of $E$ that fix $O$. By \cite[Corollary 9.4, Chapter III]{Silverman2009}, the group $\End(E)$ can be either $\ZZ$, an order in a imaginary quadratic field or an order in a quaternion algebra. In the first case, we say that $E$ does not have \textit{complex multiplication} (CM for short). In the second case $E$ has CM and in the last case $E$ is called \textit{supersingular}.

\subsection{Isogenies}

Let now $n\ge 1$ be an integer. As usual, we write $[n]\in \End(E)$ for the multiplication by $n$ map in $E$. The subgroup of $n$ torsion in $E$ is denoted by $E[n]$. Note that $E[n]=\ker([n])$. We say that an isogeny $\phi\in \End(E)$ is \textit{separable} if $\abs{\ker(\phi)} = \deg(\phi)$.

\begin{theorem}[{see {\cite[Corollary 4.11, Chapter III]{Silverman2009}}}]
\label{theorem: factorization of isogenies}
Let $E/K$ be an elliptic curve. If $\phi_1, \phi_3 \in \End(E)$, $\phi_1$ is separable and $\ker(\phi_1) \leq \ker(\phi_3)$ then there exists a unique $\phi_2 \in \End(E)$ such that $\phi_3 = \phi_2 \circ \phi_1$.
\end{theorem}

\begin{proposition}
\label{proposition: kernel of a isogeny}
Let $E/K$ be an elliptic curve and $\phi \in \End(E)$ separable. Then there exist integers $d_1,d_2 \ge 0$ with $d_1 \mid d_2$ such that 
\begin{align*}
\ker(\phi) \cong (\ZZ/d_1\ZZ) \times (\ZZ/d_2\ZZ).
\end{align*}
\end{proposition}

\begin{proof}
Let $n$ be the largest integer such that $E[n]\le \ker(\phi)$. Then by \cref{theorem: factorization of isogenies}, there is a unique $\phi'$ such that $\phi = \phi' \circ [n]$. We claim $\phi'$ has cyclic kernel. If not, there exists a prime $\ell$ such that $H = (\ZZ/\ell\ZZ)^2$ is a subgroup of $\ker(\phi')$. As every non-trivial element in $H$ has order $\ell$, we have $H \le E[\ell]$. This implies that $\ell \neq \car(K)$ by \cite[Corollary 6.4(c), Chapter III]{Silverman2009}. Then $H$ and $E[\ell]$ have the same cardinality and $H = E[\ell] \le \ker(\phi')$, so $\phi' = \phi'' \circ [\ell]$ by \cref{theorem: factorization of isogenies}. Hence $\phi=\phi''\circ [\ell n]$, which contradicts the choice of $n$. 

Since $\phi$ is separable, both $\phi'$ and $[n]$ are separable. Thus $\car(K) \nmid n$. Therefore $\ker(\phi')$ is isomorphic to $\ZZ/\deg(\phi')\ZZ$ by \cite[Theorem 4.10, Chapter III]{Silverman2009}. Moreover $E[n] \leq \ker(\phi) \leq E[n \cdot \deg(\phi')]$. This forces $\ker(\phi)$ to be a finite $2$-generated abelian group, and the result follows from the structure theorem of finite abelian groups.
\end{proof}

Let $E/K$ be an elliptic curve such that $\End(E)$ is abelian and $a \in \End(E)$. Then $\ker(a)$ is $\End(E)$-invariant as if $P \in \ker(a)$ and $\phi \in \End(E)$ then $a(\phi(P)) = \phi(a(P)) = O$. Thus, we obtain a well-defined action of $\End(E)$ on $\ker(a)$.

The following lemma gives more concrete information about this action for elliptic curves with CM: 

\begin{lemma}
\label{theorem: equivalence of End(E)-sets}
Let $E/K$ be an elliptic curve with CM and $a \in \End(E)$ separable. Then $\End(E)/a \End(E)$ is isomorphic to $\ker(a)$ as $\End(E)$-sets.
\end{lemma}
\begin{proof}
We claim that there exists $P \in \ker(a)$ such that the evaluation map
\begin{align*}
\ev_P: \End(E) \rightarrow \ker(a) \quad \text{given by}\quad \phi \mapsto \phi(P)
\end{align*}
is surjective. We prove this claim at the end of the proof. Assuming the claim, note that we have $a \End(E) \leq \ker(\ev_P)$. Therefore
\begin{align*}
\deg(a) = \abs{\ker(a)} = [\End(E): \ker(\ev_P)] \leq [\End(E): a \End(E)] = \deg(a),
\end{align*}
where in the first equality we used that $a$ is separable and in the last equality the equivalence from \cite[Section 9, Chapter III]{Silverman2009}. The equivalence as $\End(E)$-sets follows now from
\begin{align*}
\psi \circ \ev_P(\phi) = \psi(\phi(P)) = \ev_P(\psi \circ \phi)
\end{align*}
for any $\psi,\phi\in\End(E)$.

We now prove our claim. Let $\omega \in \End(E)$ be such that $\End(E) = \<[1],\omega>$ as a $\ZZ$-module. It is enough to find $P \in \ker(a)$ such that for any prime number $\ell$, the quotient $\ker(a)/\ell\ker(a)$ is generated by $\set{P, \omega(P)}\subset \IIm(\ev_P)$. Indeed, as $\ker(a)$ is a 2-generated abelian group by \cref{proposition: kernel of a isogeny}, the same holds for $\mathrm{coker}(\ev_P)=\ker(a)/\IIm(\ev_P)$. Therefore if the cokernel is non-trivial, there must exist some cyclic quotient $\pi_\ell:\ker(a)\to Q_\ell$ with $Q_\ell$ of prime order $\ell$ and with $\IIm(\ev_P)\le \ker(\pi_\ell)$. However, if $\set{P, \omega(P)}$ generates $\ker(a)/\ell\ker(a)$, then  $\pi_\ell(\IIm(\ev_P))= Q_\ell$, as $Q_\ell$ is a quotient of $\ker(a)/\ell\ker(a)$.

By \cref{proposition: kernel of a isogeny}, we have $\ker(a)\cong (\mathbb{Z}/{d_1}\mathbb{Z})\times (\mathbb{Z}/{d_2}\mathbb{Z})$ with $d_1\mid d_2$. Thus we only need to consider $\ell\mid d_2$. We will first, for each prime $\ell$, construct a point $P_\ell \in \ker(a)$ such that $\ker(a)/\ell\ker(a)$ is generated by $\{P_\ell,\omega(P_\ell)\}$. By the Chinese remainder theorem, we then find a global $P$.

If $\ell \mid d_2$ but not $d_1$, then $\ker(a)/\ell\ker(a) \cong \FF_\ell$, and we take $P_\ell$ any element outside $\ell \ker(a)$. 

If $\ell \mid d_1$, then $V_\ell:=\ker(a)/\ell\ker(a)$ is an $\mathbb{F}_p$-vector space of dimension 2. Let us write $\{e_1,e_2\}$ for a Jordan basis of $V_\ell$ with respect to the linear map $\overline{\omega}$ induced by~$\omega$. Let us write $\lambda_1,\lambda_2$ for the eigenvalues of $\overline{\omega}$.

If $\overline{\omega}$ is not diagonalizable, then $\overline{\omega}(e_2) = e_1 + \lambda_2 e_2$ and $V_\ell = \<e_2, \overline{w}(e_2)>$. Thus, choose $P_\ell$ any point in $\ker(a)$ projecting to $e_2$.

If $\overline{\omega}$ is diagonalizable and $\lambda_1\ne\lambda_2$, then $
V_\ell=\langle e_1+e_2,\overline{w}(e_1+e_2)\rangle$, as the determinant of the matrix whose columns are $\set{e_1+e_2, \overline{\omega}(e_1,e_2)}$ is $\lambda_2 - \lambda_1 \neq 0$. Thus, choose $P_\ell$ any point in $\ker(a)$ projecting to $e_1+e_2$.

To conclude, let us show that it cannot happen that $\overline{\omega}$ is diagonalizable with $\lambda_1=\lambda_2$. Choose a lift $\lambda\in \mathbb{Z}$ of the unique eigenvalue of $\overline{\omega}$ in $\mathbb{F}_\ell$. As $\ell \mid d_1$, then $V_\ell$ is isomorphic to the $\ell$-torsion of $E$. Since $\omega-[\lambda]$ maps everything to 0 in $V_\ell$, we obtain that $(\omega-[\lambda])(E([\ell]))=0$, or in other words $\ker([\ell])\le \ker(\omega-[\lambda])$. Then by \cref{theorem: factorization of isogenies}, there exists $\phi\in \End(E)$ such that $\omega-[\lambda] = \phi\circ [\ell]$. Decomposing $\phi=[m_1]+[m_2]\omega$ and comparing the coefficients of $\omega$, we get $m_2\ell=1$, which contradicts $\ell$ being prime.
\end{proof}

\subsection{Lattès maps}

Let $E/K$ be an elliptic curve and $F$ an endomorphism of $E$, i.e. $F$ does not necessarily fix $O$. By \cite[Example 4.7, Chapter III]{Silverman2009}, there exists $a \in \End(E)$ and $Q \in E(K)$ such that $F(P) = a(P) + Q$. 

Recall that $\Aut(E)$ is the group of units of the ring $\End(E)$ and let $H \leq \Aut(E)$ be a non-trivial subgroup. Let $E/H$ be the quotient curve of $E$ under the action of $H$ and $\pi: E \rightarrow E/H$ the resulting projection. The quotient map $\pi$ is always separable; see \cite[Theorem 2.4, chapter II]{Silverman2009}. By the Riemann-Hurwitz theorem, we get that the genus $g(E/H) = 0$ and therefore $E/H$ is isomorphic to $\PP^1$. 

A \textit{Lattès map} is the rational function $f \in K(z)$ such that $\pi F = f \pi$. Such a map~$f$ is well-defined if and only if $Q \in E(K)^H$. Note that $f$ is separable if and only if $F$ is separable.

\section{FPP of geometric iterated Galois groups of a Lattès maps}
\label{section: general formula}

In this section, we recall some properties of the iterated Galois groups of Lattès maps and show that their natural action is equivalent to their action on the Tate module. We conclude the section by proving a formula for the fixed-point proportion and showing that for every separable Lattès maps, its fixed-point proportion is positive.

\subsection{Fixed-point proportion}

Let $T$ be a regular rooted tree of degree $d$, i.e. a rooted tree whose vertices have exactly $d$ immediate descendants. We write $\mathcal{L}_n$ for the vertices at distance $n\ge 1$ from the root. It is well-known that $\Aut(T)$, the group of all automorphisms of $T$, is profinite. Therefore, if $G$ is a closed subgroup of $\Aut(T)$, then $G$ is also profinite and it has a unique Haar probability measure. The fixed-point proportion of $G$, denoted $\FPP(G)$, is the probability that a random element in $G$ fixes at least one vertex in $\mathcal{L}_n$ for every $n \geq 1$. If we write $\pi_n(G)$ for the induced action of $G$ on $\mathcal{L}_n$, we obtain the equivalent definition
$$\mathrm{FPP}(G):=\lim_{n\to \infty}\frac{\{g\in \pi_n(G): g\text{ fixes a vertex }v\in\mathcal{L}_n\}}{\abs{\pi_n(G)}}.$$

\begin{remark}
\label{remark: FPP infinitely many ns}
In particular, the fixed-point proportion can be computed using any subsequence of natural numbers that goes to infinity.
\end{remark}

\subsection{Structure of the geometric iterated Galois group}

Recall that $G_\infty(K,f,t)$ is the \textit{geometric iterated Galois group of $f$}, i.e. the Galois group of the infinite extension $K(\bigcup_{n\ge 1}f^{-n}(t))/K(t)$. This group acts naturally on the regular rooted tree $T$ of degree $\mathrm{deg}(f)$, obtained by taking preimages of $t$ by the iterates of $f$. We can consider the finite quotients $G_n(K,f,t)$ corresponding to the finite extensions $K(f^{-n}(t))/K(t)$.

The following theorem is well-known over the complex numbers (see for instance \cite{SelfSimilar} or \cite{Radi2026inverseGalois}) and it will appear in this generality in an upcoming paper of the second author \cite{Radi2026IGGpullbacks}:

\begin{theorem}
\label{theorem: reduction to Gn}
Let $f\in K(z)$ be a separable Lattès map. Then, there exist $m,n_0 \in \NN$ such that the subgroup $H$ commute with $a^m$ and for all $n \geq n_0$ 
\begin{align*}
G_n(K,f^m,t) \cong G_n:= \set{g: P \mapsto h(P) + k\mid  h \in H \text{ and }  k \in \ker(a^{mn})}.    
\end{align*}
Moreover, the action on the $n$th level of the tree of $G_n(K,f^m,t)$ corresponds to the action of $G_n$ on $\ker(a^{mn})$.
\end{theorem}

For most elliptic curves, the endomorphism $a$ commutes with the automorphisms in $H$. However, one can find counterexamples when $E$ is supersingular and $j(E) \in \set{0,1728}$. In view of \cref{remark: FPP infinitely many ns} and \cref{theorem: reduction to Gn}, we may replace $f$ by $f^m$ and assume from now and on that $a$ commutes with $H$. Then, $H$ acts on $\ker(a^n)$ naturally by $P \mapsto h(P)$ for every $h\in H$.

Given $h\in H$, we define the endomorphism
$$m_h:=\id-h.$$
Note that $m_h$ also acts on $\ker(a^n)$ for every $n \geq 1$.

\begin{lemma}
\label{lemma: fixed points in terms of the cosets mLambda}
Let $f\in K(z)$ be a separable Lattès map. For any $n\ge 1$, an element $g \in G_n$ given by $g(P) = h(P)+ k$ fixes a point in $\ker(a^n)$ if and only if $k \in m_h(\ker(a^n))$.
\end{lemma}
\begin{proof}
If $x_0 \in \ker(a^n)$, then $g(x_0) = x_0$ if and only if $k = (\id-h)x_0 = m_h(x_0)$.
\end{proof}

\subsection{The Tate module}

Note that $a$ maps $\ker(a^{n+1})$ to $\ker(a^n)$ for every $n\ge 1$. Thus, we may consider the inverse limit 
$$\Lambda := \varprojlim \ker(a^n).$$
This abelian group is called the \textit{Tate module} of $E$ with respect to $a$. If $\varphi \in \End(E)$ and $\varphi$ commutes with $a$, then $\varphi$ acts on each $\ker(a^n)$ inducing a map $\widehat{\varphi}: \Lambda \rightarrow \Lambda$. It is straightforward to see that if $C(a) \leq \End(E)$ is the centralizer of $a$ in $\mathrm{End}(E)$ with respect to the composition, then $\widehat{\cdot} : C(a) \rightarrow \End(\Lambda)$ is a ring homomorphism.

For $\varphi = a^n$, we obtain $\IIm(\widehat{a^n}) \supseteq \IIm(\widehat{a^{n+1}})$ and $$\Lambda/\widehat{a^n} \Lambda \cong \ker(a^n)$$
for every $n\ge 1$. Let us fix some $h\in H$. As $m_h\in C(a)$, 
\begin{align*}
|\Lambda: \IIm(\widehat{m_h})| = \lim_{n \rightarrow \infty} [\ker(a^n): \IIm(m_h|_{\ker(a^n)})] = \lim_{n \rightarrow \infty} \abs{\ker(m_h|_{\ker(a^n)})}.
\end{align*}
For any $n\ge 1$, the last term in the equation above divides $\abs{\ker(m_h)}$ which divides $\deg(m_h)$. Hence $|\Lambda: \IIm(\widehat{m_h})|$ divides $\deg(m_h)$.

For each $n \geq 1$, we define 
\begin{align*}
A_n^h := (\IIm(\widehat{a^n}) + \IIm(\widehat{m_h}))/\IIm(\widehat{m_h}).
\end{align*}
Then $A_{n+1}^h \leq A_n^h$. Define $A_\infty^h := \bigcap_{n \geq 0} A_n^h$. Notice that $A_\infty^h$ is a finite subgroup in $\Lambda/\IIm(\widehat{m_h})$. 

\subsection{A formula for the fixed-point proportion}

Now, we are in position to prove a formula for the fixed-point proportion of the geometric iterated Galois group of a Lattès map:

\begin{theorem}
\label{theorem: formula for FPP Lattès}
Let $f\in K(z)$ be a separable Lattès map. Then 
\begin{align*}
\FPP(G_\infty(K,f,t)) = \frac{1}{\abs{H}} \sum_{h \in H \setminus \set{1}} \frac{1}{|\Lambda/\IIm(\widehat{m_h}):A_\infty^h|}.
\end{align*}
In particular the fixed-point proportion is always positive.
\end{theorem}

\begin{proof}
Let $n \geq 1$ and denote by $f_n(S)$ the amount of elements in $S \subseteq G_n$ that fix at least one point in $\ker(a^n)$. Then
\begin{align*}
\FPP(G_\infty(K,f,t)) = \lim_{n \rightarrow \infty} \frac{f_n(G_n)}{\abs{G_n}} = \frac{1}{\abs{H}} \sum_{h \in H} \lim_{n \rightarrow \infty} \frac{f_n(\ker(a^n)h)}{\abs{\ker(a^n)}}.
\end{align*}
By \cref{lemma: fixed points in terms of the cosets mLambda}, we have $f_n(\ker(a^n)h) = \abs{m_h(\ker(a^n))}$, and thus 
\begin{align*}
\FPP(G_\infty(K,f,t)) = \frac{1}{\abs{H}} \sum_{h \in H} \lim_{n \rightarrow \infty} \frac{1}{|\ker(a^n):m_h(\ker(a^n))|}.
\end{align*}
If $h = 1$ then $\abs{m_h(\ker(a^n))} = 1$ and the proportion of elements in the coset $\ker(a^n)$ tends to zero when $n$ goes to infinity, so we may consider only $h \in H \setminus \set{1}$. By the definition of $\Lambda$ and the isomorphism theorems we have
\begin{align*}
\frac{\ker(a^n)}{m_h(\ker(a^n))} \cong \frac{\Lambda/(\widehat{a^n}\Lambda)}{(\IIm(\widehat{a^n})+ \IIm(\widehat{m_h}))/(\widehat{a^n}\Lambda)} \cong \frac{\Lambda}{\IIm(\widehat{a^n})+\IIm(\widehat{m_h})}\cong \frac{\Lambda/\IIm(\widehat{m_h})}{A_n^h}.
\end{align*}
Therefore
$$\lim_{n\to\infty}\frac{1}{|\ker(a^n):m_h(\ker(a^n))|}=\frac{1}{|\Lambda/\IIm(\widehat{m_h}):A_\infty^h|}$$
and the formula in the statement follows. The fixed-point proportion is always positive because
\begin{equation*}
\FPP(G_\infty(K,f,t)) \geq \frac{1}{\abs{H}} \sum_{h \in H \setminus \set{1}}\frac{1}{|\Lambda/\IIm(\widehat{m_h})|} > 0. \qedhere
\end{equation*}
\end{proof}

\section{Classification of the FPP in terms of the ramification portrait}
\label{section: classification}

In this section, we explicitly compute the fixed-point proportion of the associated iterated Galois groups in terms of the ramification portrait of each separable tamely ramified Lattés map.

\subsection{Background}

Let $e_\phi(P)$ denote the ramification degree of a point $P$ with respect to a $K$-morphism $\phi$. We say that $\phi$ is unramified at $P$ if $e_\phi(P) = 1$. We say that $\phi$ is unramified if $e_\phi(P) = 1$ for all $P$ in the domain of $\phi$. We say $\phi$ is \textit{tamely ramified} if $\car(K) \nmid e_\phi(P)$ for all $P$ in the domain of $\phi$.

Recall that if $F$ is a separable morphism of an elliptic curve $E/K$ then $F$ is unramified \cite[Theorem 4.10, Chapter III]{Silverman2009}. As $\pi$ corresponds to the morphism of a quotient curve by a non-trivial subgroup $H$, we have 
\begin{equation}
\label{equation: ramification degree quotient curve}
e_\pi(P) = \abs{\Stab_H(P)};
\end{equation}
see for instance \cite[Exercise 3.13]{Silverman2009}. In particular $e_\pi(O) = \abs{H}$. Therefore, if $f$ is a Lattès map induced by $E/K$, $H$ and $F$, then $f$ is tamely ramified if $\car(K) \nmid \abs{H}$. 

\begin{theorem}
\label{theorem: End(E) in terms of H}
Let $f$ be a separable tamely ramified Lattès map induced by an elliptic curve $E/K$, a non-trivial subgroup $H \leq \Aut(E)$, $\pi: E \rightarrow E/H$ and a morphism $F$ of $E$. Then $H$ is cyclic of order $2,3,4$ or $6$. Moreover if $E$ is not supersingular then:
\begin{enumerate}[\normalfont(i)]
\item If $H \cong C_3$ or $C_6$ then $\End(E) \cong \mathcal{O}_L$ with $L = \QQ(\sqrt{-3})$.
\item If $H \cong C_4$ then $\End(E) \cong \mathcal{O}_L$ with $L = \QQ(i)$.
\end{enumerate}
\end{theorem}

\begin{proof}
If $\Aut(E)$ is not cyclic then $\Aut(E)$ can be $\SL_2(\FF_3)$ if $\car(K) = 2$ or $\mathrm{Dic}(3)$ if $\car(K) = 3$, where $\mathrm{Dic}(3)$ is the dicyclic group; see for instance \cite[Exercise~A.1]{Silverman2009}. As $f$ is tamely ramified, if $\car(K) = 2$ then $H$ must be a subgroup of $\SL_2(\FF_3)$ of order coprime to $2$ and if $\car(K) = 3$, then $H$ must be a subgroup of $\mathrm{Dic}(3)$ of order coprime to $3$. In either case $H$ is cyclic. 

If $E$ is not supersingular and $\abs{H} > 2$ then $\End(E)$ has to be an order in an imaginary quadratic field containing $H$. As $H$ embeds into $\mathcal{O}_L$ as non-integral $|H|$-th roots of unity, the result follows.
\end{proof}

For any $f\in K(z)$, we write $C_f$ for the set of critical points of $f$ and $P_f$ for the set of postcritical points, i.e. $P_f=\bigcup_{n\ge 1}f^n(C_f)$. For $f \in K(z)$ a Lattès map, the set $P_f$ coincides with the critical values of $\pi$; see \cite[Proposition 6.45]{Silverman2007}. In particular $P_f$ is finite. 

We define the \textit{Thurston orbifold} of a Lattès map $f\in K(x)$ as the pair $(\mathbb{P}^1,\nu_f)$ for the map $\nu_f:\mathbb{P}^1\to \mathbb{N}\cup \{\infty\}$, where
$$\nu_f(z):=\mathrm{lcm}\{\nu_f(w)e_f(w) : w\in f^{-n}(z) \text{ for some }n\ge 1\}.$$
Note that $\nu_f(z)\ge 2$ only for $z\in P_f$. As $P_f$ is finite, the tuple $(\nu_f(z))_{z\in P_f}$ is finite and we call it the \textit{type} of the orbifold. We define the \textit{Euler characteristic} $\chi(\nu_f)$ of the orbifold $(\mathbb{P}^1,\nu_f)$ as
$$\chi(\nu_f):=2-\sum_{p\in P_f}\left(1-\frac{1}{\nu_f(p)}\right).$$
We say that the orbifold $(\mathbb{P}^1,\nu_f)$ is \textit{euclidean} if $\chi(\nu_f)=0$.

Lattès maps have an euclidean orbifold; see for instance \cite[Theorem 4.1]{Milnor}. Branched coverings with a euclidean orbifold were completely characterized by Douady and Hubbard in \cite{DouadyHubbard}. These two results are proved over the complex numbers but the same holds over $K$ (this will appear in the upcoming preprint \cite{SantiOrbifold}). We summarize these facts in the following theorem:

\begin{theorem}
\label{theorem: orbifolds in terms of H}
Let $f$ be a separable tamely ramified Lattès map induced by an elliptic curve $E/K$, a non-trivial subgroup $H \leq \Aut(E)$, $\pi: E \rightarrow E/H$ and a morphism $F$ of $E$. If $p \in P_f$ and $\pi(P) = p$, then $v_f(p) = e_\pi(P)$. Therefore:
\begin{enumerate}[\normalfont(i)]
\item if $H \cong C_2$ then $f$ has an euclidean orbifold of type $(2,2,2,2)$;
\item if $H \cong C_3$ then $f$ has an euclidean orbifold of type $(3,3,3)$;
\item if $H \cong C_4$ then $f$ has an euclidean orbifold of type $(2,4,4)$;
\item if $H \cong C_6$ then $f$ has an euclidean orbifold of type $(2,3,6)$.
\end{enumerate}
\end{theorem}

The next well-known result shows that the set of postcritical points of a Lattès map is critically exceptional in the sense of \cite{JorgeSanti2}:

\begin{lemma}
\label{lemma: Pf is critically exceptional}
Let $f$ be a separable tamely ramified Lattès map. Then 
$$f^{-1}(P_f) \subseteq P_f \cup C_f.$$
\end{lemma}
\begin{proof}
Recall that Lattès maps have an euclidean orbifold. Now, having an euclidean orbifold is equivalent to $\nu_f(f(z))=\nu_f(z)e_f(z)$ for every $z\in \mathbb{P}^1$; see for instance \cite[Proposition 9.1]{DouadyHubbard}. However, if there was $z\in f^{-1}(P_f)\setminus( P_f\cup C_f)$, we would reach a contradiction as
$$1<\nu_f(f(z))=\nu_f(z)e_f(z)=1.$$
Thus $f^{-1}(P_f) \subseteq P_f \cup C_f$.
\end{proof}

Given an integer $d \ge 2$, we write $\ZZ_d := \varprojlim \ZZ/d^n\ZZ$ for the ring of $d$-adic numbers. If $a$ is separable, by \cref{proposition: kernel of a isogeny} we have
\begin{align*}
\Lambda \cong \ZZ_{d_1} \times \ZZ_{d_2},
\end{align*}
for some positive integers $d_1,d_2\ge 1$ such that $d_1 \mid d_2$.

To calculate the fixed-point proportion, we need to understand the action of $\widehat{m_h}$ on $\Lambda$. If $m_h \in \ZZ$ then it acts on $\ker(a^n)$ as the scalar matrix $m_h I$ for every $n \geq 1$. This naturally lifts to the inverse limit, concluding that $\widehat{m_h}$ acts on $\Lambda$ as the same scalar matrix. If $m_h \notin \ZZ$, we will assume that $E$ is not supersingular. Then $E$ has CM and $\End(E)$ is isomorphic to one of the rings $\mathcal{O}_L$ described in \cref{theorem: End(E) in terms of H}. In particular, we shall also see $H$ as a subgroup of $\mathcal{O}_L$ in this case. By \cref{theorem: equivalence of End(E)-sets}, the groups $\ker(a^n)$ and $\End(E)/a^n\End(E)$ are isomorphic as $\End(E)$-sets, which implies that $\Lambda$ and $\varprojlim \End(E)/a^n \End(E)$ are isomorphic as $\End(E)$-sets. Identifying $\mathrm{End}(E)$ with $\mathcal{O}_L$ we get 
\begin{align*}
(\Lambda/\IIm(\widehat{m_h}))/A_n^h = \Lambda/(\IIm(\widehat{a^n})+\IIm(\widehat{m_h})) \cong \mathcal{O}_L/(m_h\mathcal{O}_L + a^n\mathcal{O}_L)
\end{align*}
for all $n \geq 1$.

When $m_h=2$, then $\mathcal{O}_L/2\mathcal{O}_L\cong E[2] \cong \mathbb{F}_2^2$ and the above quotients are no more than the cokernel of $a^n$ viewed as a linear map acting on $E[2]$. 

If the Lattès map is defined over the complex numbers, the action on the Tate module is directly given by the action on the corresponding lattice $\Lambda\subset \mathbb{C}$ which is precisely $\mathcal{O}_L$ if $E/\mathbb{C}$ has CM. In this case, we do not need to use \cref{theorem: equivalence of End(E)-sets}.

\subsection{Orbifold type $(2,2,2,2)$}

Here $H = \set{\pm 1}$ and therefore the critical points of~$\pi$ correspond to the 2-torsion $E[2]$. Since $f$ is tamely ramified and $|H|=2$, then $\mathrm{char}(K)\ne 2$ and so $E[2]\cong \mathbb{F}_2^2$; see \cite[Corollary 6.4, Chapter III]{Silverman2009}. Moreover, as $\pi F = f \pi$, the postcritical dynamics of the Lattès map are given by the action of $F$ on $E[2]$. As $Q \in E[2]$ is just a translation of all the points in $E[2]$, the ramification portrait will only depend on the action of $a$ on $E[2]$. 

As before, we only need to consider $h = -1$. Thus $m_h = 2$, and $a$ acts on $E[2]$ as a $2$-dimensional matrix $$a=\mat{a_1 & a_2 \\ a_3 & a_4}$$ with coefficients in $\mathbb{F}_2$.

Hence, the quotient $\Lambda/\IIm(\widehat{m_h})$ is precisely the cokernel of $a^n$ in $E[2]$ for large enough $n\ge 1$.

\subsubsection{First case: $\#A_h^\infty=4$}

This case corresponds to $\mathrm{coker}(a^n)=1$ for all $n\ge 1$, i.e. $a$ is invertible when viewed as a matrix with coefficients in $\mathbb{F}_2$. This is precisely the case when $a$ has odd determinant. As $\mathrm{det}(a)= \deg(f)$, this happens precisely when the Lattès map is of odd degree.

\subsubsection{Second case: $\#A_h^\infty=2$} 
\label{subsection: Ah = 2}

In this case, without loss of generality, we may assume that there exists some $n\ge 1$ such that
$$a^n=\mat{0 & 0 \\ 0 & 1}.$$ In particular, we have two distinct eigenspaces for $a^n$, meaning that the post-critical orbits must contain two distinct periodic cycles. Therefore, this happens precisely when $d$ is even and there is two disjoint postcritical orbits.

\subsubsection{Third case: $\#A_h^\infty=1$}

In this last case, there exists some $n\ge 1$ such that
$$a^n=0.$$
Thus, arguing as above, there is a unique periodic orbit in the postcritical orbits of the Lattès map as $a^n$ maps all the 2-torsion to~O. Hence, we obtain this precisely when the degree is even and all the postcritical orbits intersect.

\subsubsection{Fixed-point proportion}

We summarize the above results as follows:

\begin{theorem}[Orbifold type (2,2,2,2)]
\label{theorem: fpp (2 2 2 2)}
Let $f\in K(z)$ be a separable tamely ramified  Lattès map with orbifold of type $(2,2,2,2)$. Then
\begin{enumerate}[\normalfont(i)]
\item if $f$ is of odd degree, then $\mathrm{FPP}((G_\infty(K,f,t)))=1/2$;
\item if $f$ is of even degree and it contains two non-intersecting postcritical orbits, then $\mathrm{FPP}((G_\infty(K,f,t)))=1/4$;
\item if $f$ is of even degree and all its postcritical orbits intersect in one, then $\mathrm{FPP}((G_\infty(K,f,t)))=1/8$.
\end{enumerate}
\end{theorem}

Note that case (ii) in \cref{theorem: fpp (2 2 2 2)} only holds for elliptic curves with CM. Indeed, if $E$ does not have CM, then the isogeny $a$ gives a scalar matrix so case (ii) cannot happen. Furthermore, if $E$ does have CM both (ii) and (iii) are actually realized by Lattès maps; see for instance \cite[Proposition 2.2]{Bridy}.

\subsection{Orbifold of type (2,4,4)}

In this case $\mathcal{O}_L = \ZZ[i]$ is an euclidean domain with norm $N(a+bi)=a^2+b^2$. As remarked above, we may work in $\mathcal{O}_L$. Then, for any $h\in H$, we have $m_h\mathcal{O}_L+a^n\mathcal{O}_L$ is an ideal of $\mathcal{O}_L$ and thus 
$$m_h\mathcal{O}_L+a^n\mathcal{O}_L=\mathrm{gcd}(m_h,a^n)\mathcal{O}_L.$$

The subgroup $H$ (viewed in $\mathcal{O}_L$) is generated by $i$. Then $m_i = 1-i$ and as $N(m_i) = 2$ then $\ker(m_i) = \set{O, P_1}$. Since $\ker(m_i) \leq E[2]$ and $E[2]$ has four elements, there exists $P_2 \in E(K)$ such that $E[2] = \<P_1,P_2>$. Therefore, by \cref{equation: ramification degree quotient curve} we obtain that $\pi$ has four critical points. Moreover, the orbifold type of $f$ forces $\pi(P_2) = \pi(P_1 + P_2)$. As in the case $(2,2,2,2)$, the action of $f$ on $P_f$ is given by the action of $a$ on $E[2]$.

\subsubsection{First coset: $h = -1$}
Note that in this case, we get $m_h=2$ and we are in the same exact situation as in the type (2,2,2,2). Thus, if $f$ is of odd degree $d$ we get $\#A_h^\infty=4$ and if $d$ is even then we get either $\#A_h^\infty=2$ or $\#A_h^\infty=1$. We claim the case $\# A_h^\infty = 2$ cannot happen. Indeed if it does, then by \cref{subsection: Ah = 2} there exists $P \in E[2]$ such that $a(P) = P$, two elements in $E[2]$ are mapped to $O$ by $a$ and the other two are mapped to $P$. Moreover, by the condition $\pi(P_2) = \pi(P_1 + P_2)$, the points $P_2$ and $P_1 + P_2$ are mapped by $a$ to the same point. Then $f$ acts on $P_f$ having a fixed point $p = \pi(P)$. As $v_f(p) < \infty$ and $p$ is fixed by $f$, the Lattès map $f$ is unramified at $p$. By \cref{lemma: Pf is critically exceptional} the only preimages of $p$ are itself and some critical points, all of which have even ramification degree. Therefore $d$ is odd. However, this leads to a contradiction as $d$ odd implies $\#A_h^\infty=4$. Therefore if $d$ is even we have $\#A_h^\infty = 1$.

\subsubsection{Second coset: $h = i$}

As $N(1-i) = 2$, then $N(\mathrm{gcd}(1-i,a^n))$ is either $1$ or $2$. Observe that it is $2$ if and only if $d$ is even. Indeed, as $1-i$ is irreducible, $N(\mathrm{gcd}(1-i,a^n))=2$ is equivalent to $a^n$ being divisible by $1-i$. In $\mathbb{Z}[i]$, being divisible by $1-i$ is equivalent to having even norm. As $N(a)=d$, this is equivalent to $d$ being even.

If $d$ is even, then $N(\mathrm{gcd}(1-i,a^n))=2$ means that $(1-i)\mathcal{O}_L+a\mathcal{O}_L$ is a nilpotent ideal in $\mathcal{O}_L/(1-i)\mathcal{O}_L$. In particular, we get $\#A_h^\infty=1$.

On the other hand, if $d$ is odd, we get $N(\mathrm{gcd}(1-i,a^n)=1$ for every $n\ge 1$. In this case $\mathrm{gcd}(1-i,a^n)$ is a unit, and thus $(1-i)\mathcal{O}_L+a^n\mathcal{O}_L=\mathcal{O}_L$ for every $n\ge 1$. Therefore, in this case we get $\#A_h^\infty=2$.

\subsubsection{Third coset: $h=-i$} 

As $N(1-i)=N(1+i)=2$, the same argument as above shows that $\#A_h^\infty=2$ precisely when $d$ is odd and $\#A_h^\infty=1$ precisely when~$d$ is even.

\subsubsection{Fixed-point proportion}

Again, we summarize the results for Lattès maps with an orbifold of type (2,4,4) in the following theorem:

\begin{theorem}[Orbifold type (2,4,4)]
\label{theorem: fpp (2,4,4)}
Let $f\in K(z)$ be a separable tamely ramified Lattès map with orbifold of type $(2,4,4)$. Then
\begin{enumerate}[\normalfont(i)]
\item if $f$ is of odd degree, then $\mathrm{FPP}((G_\infty(K,f,t)))=3/4$;
\item if $f$ is of even degree, then $\mathrm{FPP}((G_\infty(K,f,t)))=5/16$.
\end{enumerate}
\end{theorem}

\subsection{Orbifold of type (3,3,3)}

In this case $\mathcal{O}_L=\mathbb{Z}[\zeta_3]$ is an euclidean domain with norm $N(a+b\zeta_3)=a^2-ab+b^2$. Thus, as before, for any $h\in H$, we have $m_h\mathcal{O}_L+a^n\mathcal{O}_L$ is an ideal of $\mathcal{O}_L$ and thus 
$$m_h\mathcal{O}_L+a^n\mathcal{O}_L=\mathrm{gcd}(m_h,a^n)\mathcal{O}_L.$$
The subgroup $H$ is generated by $\zeta_3$. 

\subsubsection{First coset: $h = \zeta_3$}

For this first coset we have $m_h=1-\zeta_3$ and thus $N(m_h)=3$. Therefore $N(\mathrm{gcd}(1-\zeta_3,a^n))\in \{1,3\}$. Again, the element $1-\zeta_3$ is irreducible and thus arguing as in the (2,4,4) case we get that $\#A_h^\infty=3$ precisely when the degree~$d$ is not divisible by 3 and $\#A_h^\infty=1$ when $d$ is divisible by 3.

\subsubsection{Second coset: $h = \zeta_3^2$} For the second coset, also $N(1-\zeta_3^2)=3$.
Hence, arguing as in the previous case, we get that $\#A_h^\infty=3$ precisely when the degree $d$ is not divisible by 3 and $\#A_h^\infty=1$ when $d$ is divisible by 3.

\subsubsection{Fixed-point proportion}

For Lattès maps with an orbifold of type (3,3,3) we get the following:

\begin{theorem}[Orbifold type (3,3,3)]
\label{theorem: fpp (3,3,3)}
Let $f\in K(z)$ be a separable tamely ramified Lattès map with orbifold of type $(3,3,3)$. Then
\begin{enumerate}[\normalfont(i)]
\item if $f$ is of degree not divisible by 3, then $\mathrm{FPP}((G_\infty(K,f,t)))=2/3$;
\item if $f$ is of degree divisible by 3, then $\mathrm{FPP}((G_\infty(K,f,t)))=2/9$;
\end{enumerate}
\end{theorem}

\subsection{Orbifold of type (2,3,6)}

In this case $\mathcal{O}_L = \mathbb{Z}[\zeta_6]$ is an euclidean domain with norm $N(a+b\zeta_6)=a^2+ab+b^2$. Then again, for any $h\in H$, we have $m_h\mathcal{O}_L+a^n\mathcal{O}_L$ is an ideal of $\mathcal{O}_L$ and thus 
$$m_h\mathcal{O}_L+a^n\mathcal{O}_L=\mathrm{gcd}(m_h,a^n)\mathcal{O}_L.$$

Now $H$ is generated by $\zeta_6$. As $f$ is tamely ramified we have $\car(K) \neq 2,3$. Therefore we can write $E: y^2 = x^3 + Ax+B$ using a Weierstrass form. As $\End(E) \cong \mathbb{Z}[\zeta_6]$, we have $j(E) = 0$ and thus we may assume that $A = 0$ and $B \neq 0$. 

By \cite[Remark 6.27]{Silverman2007} we have $[\zeta_6^k](x,y) = (\zeta_6^{4k}x, \zeta_6^{3k}y)$, so the critical points of~$\pi$ are:
\begin{enumerate}[\normalfont(i)]
\item $E[2]\setminus\set{O}$ with ramification degree $2$;
\item $\ker(1-[\zeta_6^2])\setminus\set{O}$ with ramification degree $3$;
\item and $O$ with ramification degree $6$. 
\end{enumerate}

In particular, the only point fixed by $H$ is $O$ and thus $Q \in E(K)^H=\{O\}$. Therefore $F = a$ in this case. As $a$ commutes with every $m_h$, the endomorphism $F$ preserves $E[2]$ and $\ker(1-[\zeta_6^2])$. Hence $$P_f = \set{\pi(E[2]\setminus\set{O}), \pi(\ker(1-[\zeta_6^2]\setminus\set{O}), \pi(O)}.$$ As $a$ fixes $O$, then $f$ fixes $\pi(O)$. 

\subsubsection{First and fifth cosets: $h = \zeta_6, \zeta_6^5$}

For these values of $h$, we get $N(m_h)=1$. Hence $\#A_h^\infty=1$.

\subsubsection{Second and fourth coset: $h = \zeta_6^2, \zeta_6^4$} For these cases $N(m_h)= 3$. Then, again, we get that $\#A_h^\infty=3$ precisely when the degree $d$ is not divisible by 3 and $\#A_h^\infty=1$ when $d$ is divisible by 3.

\subsubsection{Third coset: $h = \zeta_6^3$} Now $h = -1$ and we have $m_h=2$. Thus $N(m_h) = 4$.
Hence $\#A_h^\infty=4$ precisely when the degree $d$ is odd. As $E[2] \setminus \set{O}$ is all mapped by $\pi$ to the same point in $\PP^1$, the case $\# A_h^\infty = 2$ cannot happen. Therefore, if $d$ is even $\#A_h^\infty=1$.

\subsubsection{Fixed-point proportion}

Lastly, we summarize the $(2,3,6)$ case as follows:

\begin{theorem}[Orbifold type (2,3,6)]
\label{theorem: fpp (2,3,6)}
Let $f\in K(z)$ be a separable tamely ramified Lattès map of degree $d$ with orbifold of type $(2,3,6)$. Then:
\begin{enumerate}[\normalfont(i)]
\item if $d \equiv 1,5 \pmod{6}$, then $\mathrm{FPP}((G_\infty(K,f,t)))=5/6$;
\item if $d \equiv 3 \pmod{6}$, then $\mathrm{FPP}((G_\infty(K,f,t)))=11/18$;
\item if $d \equiv 2,4 \pmod{6}$, then $\mathrm{FPP}((G_\infty(K,f,t)))=17/24$;
\item if $d \equiv 0 \pmod{6}$, then $\mathrm{FPP}((G_\infty(K,f,t)))=35/72$.
\end{enumerate}
\end{theorem}

This concludes the proof of \cref{Theorem: fpp of Lattès}.


\bibliographystyle{unsrt}

\end{document}